\documentclass{amsart}

\usepackage{latexsym,amssymb,amsthm,amsmath,amscd,graphicx}

\usepackage{hyperref}
\hypersetup{
colorlinks=true,
linkcolor=blue,
anchorcolor=blue,
citecolor=blue
}

\usepackage{amsfonts}

\theoremstyle{plain}
\newtheorem{theorem}{Theorem}[section]

\newtheorem{proposition}{Proposition}[section]

\newtheorem{corollary}{Corollary}[section]

\newtheorem{definition}{Definition}[section]

\newtheorem{example}{Example}[section]
\newtheorem{counter-example}{Counter-Example}[section]
\numberwithin{equation}{section}

\theoremstyle{remark}

 \numberwithin{equation}{section}

\newtheorem*{Theorem 5.1}{{\bf Theorem 5.1}}
\newtheorem*{Theorem A}{{\bf Theorem A}}
\newtheorem*{Theorem 5.2}{{\bf Theorem 5.2}}

 \numberwithin{equation}{section}
\def\<{\left < }
\def\>{\right >}
\def\({\left ( }
\def\){\right )}

\def\r{\eqref }

\begin{document}

\title[Convex Functions are $p$-Subharmonic Functions, $p >1$ On $\mathbb{R}^n$ with Applications] {Convex Functions are $p$-Subharmonic Functions, $p > 1$ On $\mathbb{R}^n$ with Applications}

\author{Shihshu Walter Wei}
\address{Department of Mathematics\\
University of Oklahoma\\ Norman, Oklahoma 73019-0315\\ U.S.A.}
\email{wwei@ou.edu}


\date{}

\begin{abstract}
In this paper we discuss convexity, its average principle, an extrinsic average variational method in the Calculus of Variations, an average method in Partial Differential Equations, a link of convexity to $p$-subharmonicity, subsolutions to the $p$-Laplace equation, uniqueness, existence, isometric immersions in multiple settings.  In particular, we show that a convex function on $\mathbb{R}^n$ is
a $p$-subharmonic function, for every $p > 1$, and a  $C^2$ convex
function on a Riemannian manifold is a $p$-subharmonic function
$f$, for every $p > 1\, .$ We also show that a  $C^2$ convex
function which is a submersion on a Riemannian manifold is a
$p$-subharmonic function, for every $p \ge 1\, .$ 
This result is sharp. 
As further applications, via function growth estimates in $p$-harmonic geometry, we prove that every  
$p$-balanced nonnegative $C^2$ convex function on a complete noncompact Riemannian manifold is constant for $p > 1$. In particular,
every $L^q$, nonnegative, convex function of class $C^2$ on a complete noncompact Riemannian manifold is constant for $q > p -1 > 0\, .$
\end{abstract}

\maketitle
\section{Introduction}

Convex functions are fundamental objects, tools, and concepts in
various branches of mathematics.   In fact, the notion of
convexity plays an important role in several areas of Mathematics,
such as real and complex analysis, differential geometry, nonlinear
potential theory, calculus of variations, partial differential
equations, geometric measure theory, optimal control theory,
geometric function theory, and other more. Convexity enjoys the following.

\leftline{ {\bf An {\it average} principle of convexity(resp. concavity, linearity)}(cf.[W7,(8.1)])}

\noindent
Let $f$ be a convex function $($resp. concave function, linear function $)$ .  Then
\[
\begin{aligned} f(\text{average}) &\le \text{average}\, (f)\, ,\\
\big (\operatorname{resp.}\quad  f(\text{average}) &\ge \text{average}\, (f)\, ,\\
 f(\text{average}) & = \text{average}\, (f)\, \quad \big ). 
\end{aligned}
\]
Applying the above principle, where a convex function $f = \exp\, $ and ``average" is taken over two positive numbers with respect to the sum, yields
one of the simplest inequalities, G.M. $\le$  A.M. This is a sharp isoperimetric inequality for plane rectangles that has far-reaching impacts. A dual approach from discreteness to continuity yields a sharp isoperimetric inequlality for plane curves, which is equivalent to the Sobolev inequality on $\mathbb {R}^2$ with optimal constant (cf.[W7,$\S 8$]). Isoperimetric and Sobolev inequlalities are extended to Riemannian manifolds
$M$ with sharp constants and with applications to optimal sphere theorems (cf., e.g., Wei-Zhu [WZ]).
The above {\it average} principle also leads to Jensen's inequalities for $p$-Yang-Mills energy functional and for normalized exponential Yang-Mills energy functional in Gauge Theory (cf. [W7, Theorems 10.1 and 9.1]).  \smallskip

For a down-to-earth discussion, we recall a function $f
: (a,b) \to \mathbb{R}\, ,$ on an open interval $(a, b)\subset
\mathbb{R}$ is {\it convex} if for every interval $(c,d) \subset (a,b)
$, and every linear function $h$ with
$$ f = ( \text{or}\,  \le\, )\, h \ \ \text{on} \ \ \partial [c,d]\, , \ \ \text{then} \ \ f((1-\lambda)c + \lambda d) \le  h(\lambda)\, \ \ \text{on} \ \ [c,d]\, ,$$
where $h(\lambda) = (1-\lambda)f(c) + \lambda f(d)\, .$ It is an
elegant link between geometric function theory and the theory of
differential equations. Namely, if $f$ is $C^2$, then $f$ is convex
on $(a, b)$ if and only if
$$ f\quad \operatorname{is}\,  \text{a}\, {\it subsolution}\, \operatorname{of}\quad h^{\prime \prime} = 0\, , \ \ i.e. \ \ f^{\prime \prime} \ge 0 \ \  \operatorname{on} \ \ (a, b)\, .$$  

In the differentialble context the idea in a one-dimensional open interval can be extended 
to an infinite-dimensional Hibert real space $\mathcal H$. Let $f: \mathcal H \to \mathbb R$ be a function of class $C^2$ (for simplicity), and ${d^2 f_x} (v,w)$ be the second derivative of $f$ at $x \in \mathcal H$ in the directions $v,\, w,$  then we have 
\begin{proposition} $([\operatorname{ES}])$
$$\begin{aligned} &\text{If}\quad  d^2 f_x (v,v) \ge 0\quad \text{for}\, \text{all}\quad x, v \in \mathcal H\quad \text{then}\,  f\,  \text{is}\, \text{lower}\, \text{semiconinuous}\, \text{in} \,  \mathcal H.\\
&\text{Furthermore}\, , \text{the}\, \text{subsets}\, \text{of}\, \mathcal H\, \text{which}\, \text{are}\, \text{weakly}\, \text{compact}\,  \text{are}\,  \text{precisely}\, \text{those}\\
&\text{which}\,
\text{are}\, \text{bounded}\, \text{and}\, \text{weakly}\, \text{closed}.\end{aligned}
$$\end{proposition} 

This is an abstraction of the basic work of Tonelli and Morrey on convexity properties of a variational density to insure {\it existence} in the calculus of variations.(cf. [M, ES]). \bigskip

Convexity properties also lead to the {\it uniqueness} of $p$-harmonic maps of a compact Riemannian manifold into a compact manifold with nonpositive sectional curvature, without using heat flow method (cf. [W4]). When $p=2$, this generalizes the uniqueness theorem of harmonic maps due to Hartman ([H]). 
\bigskip

{\it Observing Mathematics and Nature are beautifully interwoven, and are frequently two sides of the same coin}, (as manifested by legendary sage Lao Tzu in his book Tao Te Ching,) we proposed an extrinsic, {\it average} variational method in the Calculus of Variations (cf.[W1,W2]) as an approach to confront and resolve problems in global, nonlinear analysis, geometry and physics, by which we pioneered the study of $p$-harmonic geometry (cf.e.g.,[W4,W5]), and we have found new manifolds (cf. [W3, HW, FHLW, FHW, FHJW]).
These newly found manifolds have their strong interactions with geometry, topology, analysis, partial differential equations, calculus of variations, physics, and are briefly listed in the following table (cf. [W7], page 321 and references therein).
{\fontsize{0.01}{8}\selectfont
\begin{table}[ht]
\caption{An Extrinsic {\color{red}Average} Variational Method}\label{eqtable}
\renewcommand\arraystretch{1.2}
\noindent\[
\begin{array}{|c|c|c|c|c|}
\hline
\operatorname{Mappings}&\operatorname{Functionals}&{\color{red}\operatorname{New}\, \operatorname{manifolds}\, \operatorname{found}}&\operatorname{Geometry}&\operatorname{Topology}\\
\hline
\operatorname{harmonic}\,\operatorname{map}\, \operatorname{or}&\operatorname{energy}\, \operatorname{functional}\,E\, \operatorname{or}&{\color{red}\operatorname{SSU}\, \operatorname{manifolds}\, \operatorname{or}}&\operatorname{SU}\, \operatorname{or}&
 \pi_1=\pi_2=0\\
\Phi_{(1)}-\operatorname{harmonic}\,\operatorname{map} & E_{\Phi_{(1)}}&{\color{red}\Phi_{(1)}-\operatorname{SSU}\, \operatorname{manifolds}}&\Phi_{(1)}-\operatorname{SU}& 
\pi_1=\pi_2=0\\ \hline
p-\operatorname{harmonic}\,\operatorname{map}&E_p&{\color{red}p-\operatorname{SSU}\, \operatorname{manifolds}}&p-\operatorname{SU}& \pi_1=\cdots=\pi_{[p]}=0\\
\hline
\Phi-\operatorname{harmonic}\,\operatorname{map}\, \operatorname{or}&\Phi-\operatorname{energy}\, \operatorname{functional}\,E_{\Phi}\, \operatorname{or}&{\color{red}\Phi-\operatorname{SSU}\, \operatorname{manifolds}\, \operatorname{or}}&\Phi-\operatorname{SU}\, \operatorname{or}& \pi_1=\cdots=\pi_4=0\\
\Phi_{(2)}-\operatorname{harmonic}\,\operatorname{map}&E_{\Phi_{(2)}}&{\color{red}\Phi_{(2)}-\operatorname{SSU}\, \operatorname{manifolds}}&\Phi_{(2)}-\operatorname{SU}& \pi_1=\cdots=\pi_4=0\\
\hline
\Phi_S-\operatorname{harmonic}\,\operatorname{map}&E_{\Phi_S}&{\color{red}\Phi_S-\operatorname{SSU}\, \operatorname{manifolds}}&\Phi_S-\operatorname{SU}& \pi_1=\cdots=\pi_4=0\\
\hline
\Phi_{S,p}-\operatorname{harmonic}\,\operatorname{map}&E_{\Phi_{S,p}}&{\color{red}\Phi_{S,p}-\operatorname{SSU}\, \operatorname{manifolds}}&\Phi_{S,p}-\operatorname{SU}& \pi_1=\cdots=\pi_{[2p]}=0\\
\hline
\Phi_{(3)}-\operatorname{harmonic}\,\operatorname{map}&\Phi_{(3)}-\operatorname{energy}\, \operatorname{functional}\,E_{\Phi_{(3)}}&{\color{red}\Phi_{(3)}-\operatorname{SSU}\, \operatorname{manifolds}}&\Phi_{(3)}-\operatorname{SU}& \pi_1=\cdots=\pi_{6}=0\\
\hline
\end{array}
\]
\end{table}
}
\smallskip

This extrinsic average variational method in the Calculus of Variations is
in contrast to an {\color{red} average} method in {\it Partial Differential Equations} that we applied (cf.[CW,Proposition 2.1]) to obtain sharp growth estimates for warping functions in multiply warped product manifolds, and to
solve their {\it isometric immersion problems} into Riemannian manifolds, complex space forms, quaternionic space forms, etc. ([CW]).
\bigskip

On the other hand, the idea in calculus is naturally
extended to Riemannian geometry: A smooth function $f: M \to
\mathbb{R}$ on a Riemannian manifold $M$ is said to be {\it
convex} if for each geodesic curve $c: (-\epsilon, \epsilon)\to
M\, , f \circ c : (-\epsilon, \epsilon)\to \mathbb{R}$ is a convex
function, or equivalently, its Hession $\nabla df$ on $M$ is
nonnegative definite, where $df$ is the differential of $f\, ,$
$\nabla$ is the Riemannian connection on $M\, ,$ and $\nabla
df(X,Y) \equiv (\nabla _{X} df)(Y)\, ,$ for any smooth vector
fields $X$ and $Y$ on $M$.
This idea can be extended to nonlinear potential theory: An upper
semicontinuous function $f$ is said to be {\it $p$-subharmonic}, if $f:
M \to \mathbb{R} \cup \{\infty\}\, ,$ $f \not \equiv \infty$ on
each component of $M\, ,$ and for every bounded domain $\Omega$ in
$M$, and every $p$-harmonic function $h\in \overline{\Omega}$,
with
$$f \le h
\ \ on \ \ \partial \Omega\, , \ \ then \ \ f \le h \ \ in \ \
\Omega\, .$$ Here a {\it $p$-harmonic function} $h \in
H^{1,p}_{loc}(\Omega)$ is a continuous weak solution of
$p$-Laplace equation
$$\operatorname{div}(|\nabla h |^{p-2}\nabla h ) = 0\, .$$ Again the
interplay between geometric function theory and the theory of
partial differential equations indicates that if $f$ is continuous
and $1 < p < \infty$, then $f$ is a $p$-subharmonic function in
$M$ if and only if $f$ is a {\it subsolution of the $p$-Laplace
equation}, i.e.
\begin{equation}\label{1.0} \operatorname{div} (|df|^{p-2} df ) \equiv \sum_{i=1}^n
(\nabla (|df|^{p-2} df )) (e_i, e_i) \ge 0
\end{equation} weakly in $M\, ,$ where $\{e_i\}_{i=1}^n$ is a local orthonormal frame field on $M\,
.$   That is, \begin{equation}\label{1.01} \int _M \<|df|^{p-2} df
, d\phi\>_M dv \le 0
\end{equation}
whenever $\phi \in C_0^{\infty}(M)$ is nonnegative (cf. e.g.,
[GLM],[HK],[HKM],[WLW2]). Here $\<\, , \, \>_M$ and $dv$ denote
the Rimannian metric and volume element on $M$ respectively. 
\smallskip

In fact, there are many

\begin{example}[of $p$-subharmonic functions]
Most commonly seen functions such as the
exponential function $e^x$ on $\mathbb{R}$, $e^{ |x|}$ on
$\mathbb{R}^n\, ,$ and the distance function squared in
Cartan-Hadmard manifolds are $p$-subharmonic functions for every
$p > 1\, .$ \end{example}

These examples are natural generalizations of convex
functions. In particular, we have

\begin{theorem} \label{T:1.1} A $C^2$
convex function on a Riemannian manifold $M$ is $p$-subharmonic,
for every $p
> 1\, .$
\end{theorem}

Theorem \ref{T:1.1}
is a real analog of a well-known result in complex geometry due to
Greene and Wu [GW]

\begin{theorem} \label{T:1.2} On a K\"ahler manifold, every $C^2$
convex function is plurisubharmonic.
\end{theorem}

Recall a $C^2$ real-valued function $f$ on a complex manifold is
said to be {\it plurisubharmonic} if the Levi form $Lf$ of $f$
$$Lf \equiv 4 \sum_{\alpha, \beta} \frac {\partial ^2f}{\partial z^{\alpha}\partial
\overline{z}^{\beta}}dz^{\alpha} d\overline{z}^{\beta} \ge 0\, ,$$
where $\{z^{\alpha} = x^{\alpha}+ \sqrt {-1}\,  y^{\alpha}\}$ is a
local (complex) coordinate system in $M\, ,$$dz^{\alpha} =
dx^{\alpha}+ \sqrt {-1}\,  dy^{\alpha}\, ,$ $\frac
{\partial}{\partial z^{\alpha}} = \frac 12 (\frac
{\partial}{\partial x^{\alpha}}- \sqrt {-1} \frac
{\partial}{\partial y^{\alpha}})\, ,$ and
$d\overline{z}^{\alpha}\, $ and $\frac {\partial}{\partial
\overline{z}^{\alpha}}\,  $ are complex conjugates of
$dz^{\alpha}$ and $\frac {\partial}{\partial z^{\alpha}} $
respectively.

If $f$ is a submersion, i.e. $|df| \ne 0$ everywhere, then one can
extend the range of $p\, :$

\begin{theorem} \label{T:1.3} A $C^2$
convex function on a Riemannian manifold $M$ that is a submersion,
is $p$-subharmonic  for every $p \ge 1\, .$
\end{theorem}

This result is sharp (cf. Counter-Example \ref{CE:3.1}). As
immediate consequences of Theorems \ref{T:1.1} and \ref{T:1.2}, we
have

\begin{corollary} \label{C:1.1} Every $C^2$
concave function on a Riemannian manifold $M$ is $p$-superharmonic, for any $p > 1\, ,$ and
every $C^2$ concave submersive function on $M$ is $p$-superharmonic, for
any $p \ge 1\, .$
\end{corollary}

\begin{corollary} \label{C:1.2} Let $f_i,\, i=1,2\, $ and $p$ be as in the assumption and conclusion of Theorem \ref{T:1.1} or \ref{T:1.2}
respectively.
Let $\lambda
> 0\, ,$ then $\lambda f_1 $, $f_1 + f_2\, ,$ and $\max \{f_1,
f_2\}$ are $p$-subharmonic functions.
\end{corollary}

\begin{corollary} \label{C:1.3} Let an increasing sequence of functions $\{f_i\}_{i=1}^{\infty}$ and $p$ be as in the Corollary \ref{C:1.2}.
Then $\lim_{i \to \infty} f_i $ is $p$-subharmonic.
\end{corollary}

If $M$ is Euclidean space $\mathbb{R}^n$, then one can drop the
$C^2$ assumption on $f$:

\begin{theorem} \label{T:1.4} A
convex function on $\mathbb{R}^n$ is $p$-subharmonic, for every $p
> 1\, ,$ and a convex function on $\mathbb{R}^n$ with the $n-$dimensional Lebesgue
measure $\mathcal{L}^n(\{x \in \mathbb{R}^n: |d f|=0\})=0\, ,$ is
$p$-subharmonic, for every $p \ge 1\, ,$
\end{theorem}

This result is sharp (cf. Counter-Example \ref{CE:3.1}).

In this paper we combine the link between convex functions and $p$-subharmonic functions, and the estimates on the growth of $p$-subharmonic functions (cf.[WLW2], or \S2) to prove Liouville type theorems for convex functions. 
We recall for a given $q\in \mathbb R$, a function or a differential form or a bundle-valued differential form $f$ has $p$-balanced growth (or, simply, is $p$-balanced) if $f$ has one of the following: $p$-finite, $p$-mild, $p$-obtuse, $p$-moderate, or $p$-small growth, and has $p$-imbalanced growth (or, simply, is $p$-imbalanced) otherwise (cf. [WLW2], or $\S 2$ ).
As further applications, we have the following. \medskip

\noindent {\bf Theorem 5.1} [Louiville Type Theorem for Convex Functions]  {\it Every $p$-balanced nonnegative $C^2$ convex function on a complete noncompact Riemannian manifold is constant for $p > 1\, .$ }\medskip

\noindent {\bf Corollary 5.1} {\it Every $L^q$-  nonnegative $C^2$ convex function on a complete noncompact Riemannian manifold is constant for $q > p - 1 >0\, .$  }

\section{Preliminaries}

Let $(M,g)$ be a smooth Riemannian manifold. Let $\xi :E\rightarrow M$ be a smooth Riemannian vector bundle over $(M,g)\, ,$ i.e. a vector bundle such that at each fiber is equipped with a positive inner product $\langle \quad , \quad \rangle_E\, .$
 Set $A^k(\xi )=\Gamma (\Lambda
^k T^{*}M\otimes E)$ the space of smooth $k-$forms on $M$ with
values in the vector bundle $\xi :E\rightarrow M$.  For two forms $\Omega ,\Omega ^{\prime }\in A^k(\xi )$, the
induced inner product $ \langle\Omega ,\Omega ^{\prime }\rangle$ is defined as in 
follows:
\begin{equation*}
\aligned
 \langle\Omega ,\Omega ^{\prime }\rangle&=\sum_{i_1<\cdots <i_k}\langle\Omega
(e_{i_1},...,e_{i_k}),\Omega ^{\prime }(e_{i_1},...,e_{i_k})\rangle_E \\
&=\frac 1{k!}\sum_{i_1,...,i_k}\langle\Omega
(e_{i_1},...,e_{i_k}),\Omega ^{\prime
}(e_{i_1},...,e_{i_k})\rangle_E\, ,
\endaligned\label{6.4}
\end{equation*}
where $\{e_1, \cdots e_n\}$ is a local orthonormal frame field on $(M,g)$\, . For $\Omega \in A^k(\xi )$, set  $|\Omega|^2 = \langle \Omega, \Omega \rangle$. Then $|\Omega|^q = \langle \Omega, \Omega \rangle ^{\frac q2}$. Following [W6], we introduce the following notions. 
\begin{definition}\label {D: 5.1} For a given $q\in \mathbb R\, ,$ a function or a differential form or a bundle-valued differential form $f$ has \emph{$p$-{finite growth}} $($or, simply, \emph{is
$p$-{finite}}$)$ if there exists $x_0 \in M$ such that
\begin{equation}
\liminf_{r\rightarrow\infty}\frac{1}{r^p}\int_{B(x_0;r)}|f|^{q}\, dv
<\infty\, , \label{5.1}
\end{equation}
and has \emph{$p$-{infinite growth}} $($or, simply, \emph{is
$p$-infinite}$)$ otherwise. \smallskip

For a given $q\in \mathbb R\, ,$ a  function or a differential form or a bundle-valued differential form $f$ has \emph{$p$-mild growth}
$($or, simply, \emph{is $p$-mild}$)$ if there exist  $ x_0 \in M\,
,$ and a strictly increasing sequence of $\{r_j\}^\infty_0$ going
to infinity, such that for every $l_0>0$, we have
\begin{equation}
\begin{array}{rll}
\sum\limits_{j=\ell_0}^{\infty}
\bigg(\frac{(r_{j+1}-r_j)^p}{\int_{B(x_0;r_{j+1})\backslash
B(x_0;r_{j})}|f|^q\, dv}\bigg)^{\frac1{p-1}}=\infty \, ,
\end{array}    \label{5.2}
\end{equation}
and has \emph{$p$-severe growth} $($or, simply, \emph{is
$p$-severe}$)$ otherwise. \smallskip

For a given $q\in \mathbb R\, ,$ a function or a differential form or a bundle-valued differential form $f$ has \emph{$p$-obtuse growth}
$($or, simply, \emph{is $p$-obtuse}$)$ if there exists $x_0 \in M$
such that for every $a>0$, we have
\begin{equation}
\begin{array}{rll}
\int^\infty_a\bigg( \frac{1}{\int_{\partial
B(x_0;r)}|f|^qds}\bigg)^\frac{1}{p-1}dr =
 \infty \, ,   \label{5.3}
\end{array}
\end{equation}
and has \emph{$p$-acute growth} $($or, simply, \emph{is
$p$-acute}$)$ otherwise. \smallskip

For a given $q\in \mathbb R\, ,$ a function or a differential form or a bundle-valued differential form $f$ has \emph{$p$-moderate
growth} $($or, simply, \emph{is $p$-moderate}$)$ if there exist  $
x_0 \in M$, and $\psi(r)\in {\mathcal F}$, such that
\begin{equation}\label{5.4}
\limsup _{r \to \infty}\frac {1}{r^p \psi^{p-1} (r)}\int_{B(x_0;r)}
|f|^{q}\, dv < \infty \, ,
\end{equation}
and has \emph{$p$-immoderate growth} $($or, simply, \emph{is
$p$-immoderate}$)$ otherwise, where
\begin{equation}\label{5.5} {\mathcal F} = \{\psi:[a,\infty)\longrightarrow
(0,\infty) |\int^{\infty}_{a}\frac{dr}{r\psi(r)}= \infty \ \ for \ \
some \ \ a \ge 0 \}\, .\end{equation} $($Notice that the functions
in {$\mathcal F$} are not necessarily monotone.$)$ \smallskip

For a given $q\in \mathbb R\, ,$ a function or a differential form or a bundle-valued differential form $f$ has \emph{$p$-small growth}
$($or, simply, \emph{is $p$-small}$)$ if there exists $ x_0 \in
M\, ,$ such that for every $a
>0\, ,$we have
\begin{equation}
\begin{array}{rll}
\int
_{a}^{\infty}\bigg(\frac{r}{\int_{B(x_0;r)}|f|^{q}\, dv}\bigg)^{\frac1{p-1}}
dr = \infty \, ,
\end{array}    \label{5.6}
\end{equation}
and has \emph{$p$-large growth} $($or, simply, \emph{is
$p$-large}$)$ otherwise. \end{definition}

\begin{definition}\label {D: 5.2} For a given $q\in \mathbb R\, ,$ a function or a differential form or a bundle-valued differential form $f$ has
\emph{$p$-balanced growth} $(or, simply, \emph{is $p$-balanced})$
if $f$ has one of the following: \emph{$p$-finite},
\emph{$p$-mild}, \emph{$p$-obtuse}, \emph{$p$-moderate}, or
\emph{$p$-small} growth, and has \emph{$p$-imbalanced growth} $( $
or, simply,
is $p$-\emph{imbalanced}$)$ otherwise.
\end{definition}

The above definitions of ``$p$-balanced, $p$-finite, $p$-mild,
$p$-obtuse, $p$-moderate, $p$-small" and their counter-parts
``$p$-imbalanced, $p$-infinite, $p$-severe, $p$-acute,
$p$-immoderate, $p$-large" growth depend on $q$, and $q$ will be
specified in the context in which the definition is used.

\begin{theorem} [$\operatorname{[W6]}$, $\operatorname{Theorem}$ 5.4.]
For a given $q \in \mathbb R\, ,$ a function, or differential form or bundle-valued differential form $f$ is
\[
%\begin{aligned}
 p-moderate\, \eqref{5.4}\quad  \Leftrightarrow \quad p-small\, \eqref{5.6}\quad \Rightarrow \quad p-
mild\, \eqref{5.2}\quad \Rightarrow \quad p-obtuse\, \eqref{5.3} \]
%& 
%\[ \operatorname{or}\quad \operatorname{equiavalently},\\
%& 
or equiavalently, 
\[ p-acute \quad   \Rightarrow \quad p-severe \quad \Rightarrow \quad p-large \quad \Leftrightarrow \quad p-immoderate.
%\end{aligned}
\]
Hence, for a given $q \in \mathbb R\, ,$ $f$ is 
\[
%\begin{aligned}
\quad p-\operatorname{balanced}\quad
%& 
\Rightarrow \quad \operatorname{either} \quad p-\operatorname{finite}\, \eqref{5.1}\quad \operatorname{or}\quad p-\operatorname{obtuse}\, \eqref{5.3}\]
\[
\quad p-\operatorname{imbalanced}\quad
%& 
\Rightarrow \quad \operatorname{both} \quad p-\operatorname{infinite}\quad \operatorname{and}\quad p-\operatorname{immoderate}.
%\end{aligned}
\]
If in addition, $\int_{B(x_0;r)}|f|^{q}dv$
 is convex in $r$, then the following four types of growth are all equivalent: $f$ is \emph{$p$-mild}, \emph{$p$-obtuse}, \emph{$p$-moderate}, and
\emph{$p$-small} $(\operatorname{resp.}$ 
\emph{$p$-severe}, \emph{$p$-acute}, \emph{$p$-immoderate}, and 
\emph{$p$-large}$)$, i.e.,  $f$ is
$$\qquad \eqref{5.2}\quad  \Leftrightarrow \quad \eqref{5.3}\quad \Leftrightarrow \quad \eqref{5.4}\quad \Leftrightarrow \quad \eqref{5.6}\quad $$for the same value of $q \in \mathbb R\, .$
 \label{T:2.1} 
 \end{theorem}
In particular, we have\vskip0.1in
\begin{corollary}$([W6, \operatorname{Corollary} 5.1])$\label{C:2.1}
 Every $L^q$ function or differential form or bundle-valued differential form $f$ on $M$ has \emph{$p$-balanced} growth, $p \ge 0\, ,$ and in fact, has \emph{$p$-finite},
\emph{$p$-mild}, \emph{$p$-obtuse}, \emph{$p$-moderate}, and
\emph{$p$-small} growth, $p \ge 0\, ,$  for the same value of $q$
\end{corollary}

In $[\operatorname{WLW}2]$, among many different types of inequalities on a complete noncompact Reimannian manifold $M$, we have the following uniqueness property.
\bigskip

\begin{theorem}[Liouville Property for solutions of $f \text{div}(|\nabla f|^{p-2}\nabla f)\ge 0$]
Every $C^2$ solution $ f : M \to (-\infty, \infty)$ of $f \operatorname{div}(|\nabla f|^{p-2}\nabla f)\ge 0$ is constant provided f is $p$-balanced, i.e. $f$ is one of the following:
\emph{$p$-finite}, \emph{$p$-mild}, \emph{$p$-obtuse},
\emph{$p$-moderate}, or \emph{$p$-small}, for some $q>p-1$. In
particular, every $C^2$, $L^q$ solution f of $f \text{div}(|\nabla f|^{p-2}\nabla f)\ge 0$
is constant for any $q
> p-1$. \label{T: 2.2}
\end{theorem} 
\section{Convexity and $p$-Subharmonicity }

For completeness, we prove theorems that link convexity and $p$-subharmonicity as follows. \smallskip 

\noindent {\bf Proof of Theorem \ref{T:1.1}.} Let
$\{e_i\}_{i=1}^n$ be a local orthonormal frame field on $M\, .$
Then the differential $df$ satisfies
\begin{equation}\label{1.1} \sum_{i=1}^n df(e_i)e_i = \sum_{i=1}^n
\<\nabla f, e_i\>_M e_i = \nabla f
\end{equation}
At point $x\in M$ where $|\nabla f|\ne 0\, ,$ we may assume,
without loss of generality $e_1 = \frac {\nabla f}{|\nabla f|}$.
\begin{equation}\label{1.2}
\begin{array}{rll}
(\nabla_{\nabla f}  df)( \nabla f )
& = (\nabla  df)(\nabla f, \nabla f ) \\
& = (\nabla  df)(
|\nabla f| e_1, |\nabla f| e_1)\\
& = |df|^2(\nabla df)( e_1, e_1)
\end{array}
\end{equation}
At point $x\in M$ where $|\nabla f|= 0\, ,$
\begin{equation}\label{1.3}
\begin{array}{rll}
(\nabla_{\nabla f}  df)( \nabla f ) & = 0\, .
\end{array}
\end{equation}
Let $j>0$ be an integer. It follows from \r{1.1} and \r{1.2} that

\begin{equation}\label{1.4}
\begin{array}{rll}
& \quad \operatorname{div}((|df|^2 + \frac 1j)^{\frac {p-2}2} df)\\
& \equiv \sum_{i=1}^n \nabla ((|df|^2 + \frac 1j)^{\frac {p-2}2} df) (e_i, e_i)\\
&  = \sum_{i=1}^n (\nabla _{e_i}((|df|^2 + \frac 1j)^{\frac
{p-2}2}
df)) (e_i)\\
&= \sum_{i=1}^n \bigg ( \big ({e_i} (|df|^2 + \frac 1j)^{\frac {p-2}2} \big ) df \bigg ) (e_i) + \sum_{i=1}^n ((|df|^2 + \frac 1j)^{\frac {p-2}2} \nabla _{e_i} df) (e_i)\\
& = \sum_{i=1}^n( (p-2)(|df|^2 + \frac 1j)^{\frac {p-4}2}\<\nabla _{e_i} df, df \> df) (e_i) + \sum_{i=1}^n ((|df|^2 + \frac 1j)^{\frac {p-2}2} \nabla  df)({e_i}, e_i)\\
& = \sum_{i=1}^n (p-2)(|df|^2 + \frac 1j)^{\frac {p-4}2}\<\nabla _{e_i} df, df \> df (e_i) + \sum_{i=1}^n ((|df|^2 + \frac 1j)^{\frac {p-2}2} \nabla  df)({e_i}, e_i)\\
& = (p-2)(|df|^2 + \frac 1j)^{\frac {p-4}2}\<\nabla _{\nabla f} df, df \>   + \sum_{i=1}^n ((|df|^2 + \frac 1j)^{\frac {p-2}2} \nabla  df)({e_i}, e_i)\\
& = (p-2)(|df|^2 + \frac 1j)^{\frac {p-4}2}(\nabla _{\nabla f} df)(\nabla f )   + \sum_{i=1}^n ((|df|^2 + \frac 1j)^{\frac {p-2}2} \nabla  df)({e_i}, e_i)\\
\end{array}
\end{equation}

\noindent At point $x\in M$ where $|\nabla f|\ne 0\, ,$ for every
integer $j > 0\, ,$ the last expression
\begin{equation} \label{1.5}\begin{array}{rll}
&(p-2)(|df|^2 + \frac 1j)^{\frac {p-4}2}(\nabla _{\nabla f}
df)(\nabla f )   + \sum_{i=1}^n ((|df|^2 + \frac 1j)^{\frac
{p-2}2} \nabla  df)({e_i},
e_i)\\
& = (p-2)(|df|^2 + \frac 1j)^{\frac {p-4}2}|df|^2(\nabla df)(e_1,
e_1 )   + \sum_{i=1}^n (|df|^2 + \frac 1j)^{\frac {p-4}2}(|df|^2 +
\frac 1j)( \nabla
df)({e_i}, e_i) \\
& \ge (p-1)(|df|^2 + \frac 1j)^{\frac {p-4}2}|df|^2(\nabla
df)(e_1, e_1 )   + \sum_{i=2}^n (|df|^2 + \frac 1j)^{\frac
{p-2}2}( \nabla
df)({e_i}, e_i) \\
&\ge 0\end{array}
\end{equation}
 by \r{1.2}, the convexity of $f$ and $p \ge 1$.

\noindent At point $x\in M$ where $|\nabla f| = 0\, ,$ for every
$j > 0\, ,$ the last expression in \r{1.4}
\begin{equation} \label{1.6}\begin{array}{rll}
&(p-2)(|df|^2 + \frac 1j)^{\frac {p-4}2}(\nabla _{\nabla f}
df)(\nabla f ) + \sum_{i=1}^n ((|df|^2 + \frac 1j)^{\frac {p-2}2}
\nabla df)({e_i},
e_i)\\
 &\ge 0\end{array}
\end{equation}
 by \r{1.3}, the convexity of $f$ and $p \ge 1$. Combining \r{1.4}, \r{1.5}, and \r{1.6}, we have
for every integer $j > 0\, ,$ \begin{equation} \label{1.7}\operatorname{div}
((|df|^2 + \frac 1j)^{\frac {p-2}2} df)\ge 0 \ \ everywhere \ \ in
\ \ M.\end{equation}  It follows from the Stoke's Theorem that for
every integer $j > 0\, ,$
\begin{equation}\label{1.8}
\int _M \<(|df|^2 + \frac 1j)^{\frac {p-2}2} df, d\phi\>_M dx \le
0
\end{equation}
whenever $\phi \in C_0^{\infty}(M)$ is nonnegative.

\noindent
As Cauchy-Schwarz inequality yields $$\<(|df|^2 + \frac 1j)^{\frac {p-2}2} df, d\phi\>_M \le
(|df|^2 + \frac 1j)^{\frac {p-2}2}|df| |d\phi| \le (|df|^2 +
1)^{\frac {p-1}2} |d\phi| \in L^1(M)$$ It follows from the fact
that \begin{equation}\label{1.9}\lim_{j\to \infty}(|df|^2 + \frac
1j)^{\frac {p-2}2} df = |df|^{p-2} df\, \end{equation} everywhere
for $p > 1\, ,$ the dominated convergence theorem and \r{1.8}, we
obtain the desired
\begin{equation}
\label{1.10}\int _M \<|df|^{p-2} df, d\phi\>_M dx = \lim_{j\to
\infty}\int _M \<(|df|^2 + \frac 1j)^{\frac {p-2}2} df, d\phi\>_M
dx \le 0\, .\end{equation}

\noindent {\bf Proof of Theorem \ref{T:1.3}.} If $f$ is a
submersion, then for $p \ge 1\, ,$ \r{1.9} holds  and hence
\r{1.10} completes the proof.
\smallskip

\noindent {\bf Proof of Corollary \ref{C:1.1}} Since $f$ is
concave, $-f$ is convex and hence by Theorem \ref{T:1.1}, for
every $p \ge 1\, ,$ $-f$ is a $p$-subharmonic function, or $f$ is
a $p$-superharmonic function.
\smallskip

\noindent {\bf Proof of Theorem \ref{T:1.4}.} If $f$ is a convex
function on $\mathbb{R}^n\, ,$ then by Aleksandrov's Theorem $($cf.
$[\operatorname{R}],[\operatorname{EG}])$, $f$ has a second derivative  $\mathcal{L}^n$ almost
everywhere.  Replacing $M$ in the proofs of Theorem \ref{T:1.1}
and \ref{T:1.3} with $\mathbb{R}^n$, and ``everywhere" with ``
$\mathcal{L}^n$ a.e." complete the proof.
\smallskip

\section{A Counter-Example}
In this section, we  show the optimality of $p \ge 1$ in Theorem
\ref{T:1.3} $($in which $M = \mathbb{R}^n \backslash \{0\})$ and
Theorem \ref{T:1.4} $($in which $|df| \ne 0 \ \ \mathcal{L}^n$ a.e.$)$
by giving

\begin{counter-example}\label{CE:3.1}
The function $f(x)=e^{|x|^2}$ in  Euclidean space $\mathbb{R}^{n}$
is a convex function and is not a $p$-subharmonic function for any
$p<1\, .$\end{counter-example}

\begin{proof}
By a straightforward computation, we have:
\begin{eqnarray*}
\operatorname{div}(|df|^{p-2}df) &=& \operatorname{div}\big((2re^{r^2})^{p-2}2x_{1}e^{r^2},
\cdots, (2re^{r^2})^{p-2}2x_{n}e^{r^2}\big)
\nonumber\\
&=&
2^{p-1}\sum_{i=1}^{n}\{e^{(p-1)r^2}r^{p-2}+e^{(p-1)r^2}(p-2)r^{p-3}\frac{x_i^2}{r}+e^{(p-1)r^2}(p-1)2r^{p-2}x^2_i\}
\nonumber\\
&=&
2^{p-1}\{nr^{p-2}e^{(p-1)r^2}+(p-2)r^{p-2}e^{(p-1)r^2}+2(p-1)r^{p}e^{(p-1)r^2}\}
\nonumber\\
&=& (n+p-2+2(p-1)r^2)2^{p-1}e^{(p-1)r^2}r^{p-2}\\ & <& 0
\end{eqnarray*}
for sufficiently large $r > 0\, ,$ if $\, p<1\, .$
\end{proof}
The above computation also shows that $f(x)=e^{|x|^2}$ in
$\mathbb{R}^{n}$ is a $p$-subharmonic function for every $p\ge 1\,
.$
\bigskip

\section{Further Applications}
In this section, we utilize the link between convex functions and $p$-subharmonic functions, and apply the estimates on the growth of $p$-subharmonic functions in $[\operatorname{WLW2}]$ to prove Liouville type theorem for convex functions.

\begin{theorem} [Louiville Type Theorem for Convex Functions] 
Every $p$-balanced nonnegative $C^2$ convex function on a complete noncompact Riemannian manifold $M$ is constant for $p > 1\, .$
\label{T:5.1}
\end{theorem}
\bigskip

\noindent {\bf Proof of Theorem \ref{T:5.1}.} Since $f$ is a $C^2$ convex function on $M$, Theorem \ref{T:1.1} implies that $f$ is a $p$-subharmonic function for $p > 1$. 
This is equivalent to $f$ is a subsolution of the $p$-Laplace equation, i.e., $\text{div}(|df|^{p-2})df \ge 0$. In view of $f \ge 0$, we have $f \text{div}(|df|^{p-2})df \ge 0\, .$ It follows from Theorem \ref{T: 2.2} that $f$ is constant. \bigskip

\begin{corollary}
Every $C^2$, $L^q$ convex function on a complete noncompact Riemannian manifold $M$
is constant for any $q
> p-1 > 0$
\end{corollary}

\begin{proof}
This follows from Corollary 2.1. that if $f$ is in $L^q\, ,$ then $f$ is $p$-balanced, $p > 0$ for the same $q$ (cf. $[\operatorname{W5}]$). So we can apply Theorem \ref{T:5.1}, and the result follows.
\end{proof}

\centerline{References}
\bigskip

\smallskip\noindent$[CW]$ B.-Y. Chen and S.W. Wei, {\it Sharp growth estimates for warping functions in multiply warped product manifolds}. J. Geom. Symmetry Phys. {\bf 52} (2019), 27-46.

\smallskip\noindent$[DW]$
Yuxing Dong and Shihshu Walter Wei, {\em On vanishing theorems for vector bundle valued $k$-forms and their
applications},  Communications in Mathematical Physics 304, no. 2, (2011), 329-368.
arXiv:1003.3777 .\smallskip

\smallskip\noindent{$[EG]$} L.C. Evans and R.F. Gariepy, {\em Measure theory and fine
properties of functions}, Studies in Advanced Mathematics, CRC
Press, Boca Raton, 1992

\smallskip\noindent{$[ES]$} Eells, James, Jr.; Sampson, Joseph H. Variational theory in fibre bundles. 1966 Proc. U.S.-Japan Seminar in Differential Geometry (Kyoto, 1965) pp. 22-33 Nippon Hyoronsha, Tokyo

\smallskip\noindent{$[FHLW]$}  S. Feng, Y. Han, X. Li and S.W. Wei, 
{\em The geometry of $\Phi_{S}$-harmonic maps},
J. Geom. Anal.31(2021), no.10, 
9469–9508.

\smallskip\noindent{$[FHW]$}   S. Feng, Y. Han and S.W. Wei,
{\em Liouville type theorems and stability of $\Phi_{S,p}$-harmonic maps},
Nonlinear Anal.212 (2021), Paper No. 112468, 38 pp.

\smallskip\noindent{$[FHJW]$}   S. Feng, Y. Han, K. Jiang and S.W. Wei,  
{\em The geometry of  $\Phi_{(3)}$-harmonic maps},
Nonlinear Anal. 234 (2023), Paper No. 113318, 38 pp.; arXiv:2305.19503

\smallskip\noindent{$[GLM]$} S. Granlund, P.
Lindqvist, O. Martio,{\em Conformally invariant variational
integrals}, Trans. Amer. Math. Soc. 277 (1983), no. 1, 43--73.

\smallskip\noindent{$[GW]$} R.E. Greene and H. Wu {\em On the subharmonicity and
plurisubharmonicity of geodesically convex functions}, Indiana
Univ. Math. J. 22 (1972/73), 641--653

\smallskip\noindent{$[H]$} P. Hartman, On homotopy harmonic maps, Canad. J. Math., 19 (1967), 673–687.

\smallskip\noindent{$[HK]$} J. Heinonen and T. Kilpeläinen, {\em $A$-superharmonic functions and supersolutions of degenerate elliptic equations}, Ark. Mat. 26 (1988), no. 1, 87--105

\smallskip\noindent{$[HKM]$} I. Holopainen, T. Kipel\"ainen and O. Martio, {\em Nonlinear
potential theory of degenerate elliptic equations}, Oxford
Mathematical Monographs, Clarendon Press, Oxford-New York-Tokyo,
(1993).

\smallskip\noindent{$[HW]$} Han, Yingbo; Wei, Shihshu Walter $\Phi$-harmonic maps and $\Phi$-superstrongly unstable manifolds. J. Geom. Anal. 32 (2022), no. 1, Paper No. 3, 43 pp.

\smallskip\noindent{$[LW]$} W.P. Li and S.W. Wei
Geometry and topology of submanifolds and currents.
Selected papers from the 2013 Midwest Geometry Conference (MGC XIX) held at Oklahoma State University, Stillwater, OK, October 19, 2013 and the 2012 Midwest Geometry Conference (MGC XVIII) held at the University of Oklahoma, Norman, OK, May 12-13, 2012. Edited by Weiping Li and Shihshu Walter Wei. Contemporary Mathematics, 646. American Mathematical Society, Providence, RI, 2015. ix+186 pp.\smallskip

\smallskip\noindent{$[M]$} C.B. Morrey {\it Multiple integrals in the calculus of variations}, Collog.
Lectures A.M.S. 1964.

\smallskip\noindent{$[R]$} J. G. Re\v setnjak, {\em Generalized derivatives and differentiability almost everywhere}, (Russian) Mat. Sb. (N.S.) 75(117) 1968 323--334

\smallskip\noindent{$[W1]$} S.W. Wei {\em An average
process in the calculus of variations and the stability of
harmonic maps}, Bulletin, Institute of Mathematics, Academia Sinica, 11 (1983),
no. 3, 469--474. \smallskip

\smallskip\noindent{$[W2]$} S.W. Wei {\em An extrinsic average variational method}, Recent developments
in geometry $($Los Angeles, CA, 1987$)$, edited by R. Greene, S.Y.
Cheng and H.I. Choi, 55--78, Contemporary. Math., 101, American Mathematical
Society, Providence, RI, 1989. 

\smallskip\noindent{$[W3]$} S.W. Wei {\em Liouville theorems and regularity of minimizing harmonic maps into super-strongly unstable manifolds}, Geometry and nonlinear partial differential equations (Fayetteville, AR, 1990), 131–154.
Contemp. Math., 127 American Mathematical Society, Providence, RI, 1992

\smallskip\noindent{$[W4]$} S.W. Wei, {\em Representing homotopy
groups and spaces of maps by $p$-harmonic maps}. Indiana Univ.
Math. J. 47 (1998), no. 2, 625--670. \smallskip

\smallskip\noindent{$[W5]$} S.W. Wei, {\em The unity of $p$-harmonic geometry},
Recent developments in geometry and analysis, pp.\,439-483, Advanced Lecture Math. (ALM) 23, International Press, Somerville, MA, (2012)\smallskip

\smallskip\noindent{$[W6]$} S. W. Wei, {\it Dualities in Comparison Theorems and Bundle-Valued Generalized Harmonic Forms on Noncompact Manifolds}, Sci. China Math. 64 (2021), no. 7, 1649-1702. 54 pp.

\smallskip\noindent{$[W7]$} S. W. Wei, {\it On exponential Yang-Mills fields and $p$-Yang-Mills fields}
Adv. Anal. Geom., 6, De Gruyter, Berlin, (2022) 317-358; https:$//$doi.org/10.1515/9783110741711-018; arXiv:2205.03016

\smallskip\noindent{$[WLW1]$} S. W. Wei,  J.F. Li, and L. Wu, {\em Convex functions are $p$-harmonic functions, $p > 1$ in $\mathbb R^n$}, Global Journal of Pure and Applied Math. (2007), no. 3, 219–225.

\smallskip\noindent{$[WLW2]$} S. W. Wei,  J.F. Li, and L. Wu, {\em Generalizations of the uniformization theorem and
Bochner's method in $p$-harmonic geometry}, Commun. Math. Anal. 2008, Conference 1, 46–68.

\smallskip\noindent{$[WZ]$} S.W. Wei and M. Zhu, Sharp isoperimetric inequalities and sphere theorems. Pacific J. Math. {\bf 220}(2005), 183–195.

\enddocument